\documentclass[12pt, final, reqno]{amsart}
\usepackage[mathscr]{eucal}

\theoremstyle{plain}
\newtheorem*{theorem}{Theorem}
\newtheorem{lemma}{Lemma}
\newtheorem{proposition}[lemma]{Proposition}

\theoremstyle{remark}

\newcommand{\al}{_{\alpha}}
\newcommand{\be}{_{\beta}}
\newcommand{\C}{\mathbb C}

\newcommand{\cBH}{\mathscr{B}(\mathscr{H})}
\newcommand{\cF}{\mathscr{F}}
\newcommand{\cFH}{\mathscr{F}(\mathscr{H})}
\newcommand{\cH}{\mathscr{H}}
\newcommand{\cM}{\mathscr{M}}
\newcommand{\la}{\langle}
\newcommand{\lam}{\lambda}
\newcommand{\N}{\mathbb N}
\newcommand{\od}{\odot}
\newcommand{\ot}{\otimes}
\newcommand{\cPH}{\mathscr{P}(\mathscr{H})}
\newcommand{\ra}{\rangle}
\newcommand{\rng}{\operatorname{rng}}
\newcommand{\tr}{\operatorname{tr}}

\begin{document}
\vglue -1cm
\centerline{\large{\textbf{A GENERALIZATION OF WIGNER'S
UNITARY-ANTIUNITARY}}}

\centerline{\large{\textbf{THEOREM TO HILBERT MODULES}}}
\vskip .5cm
\centerline{LAJOS MOLN\'AR}

\centerline{Institute of Mathematics and Informatics}

\centerline{Lajos Kossuth University}

\centerline{4010 Debrecen, P.O.Box 12}

\centerline{Hungary}

\centerline{phones: direct ++ 52 316 666/2815; dept. ++ 52 489 100;
                    fax ++ 52 416 857}

\centerline{e-mail: \texttt{molnarl@math.klte.hu}}

\vskip .5cm
\centerline{Running title:}

\centerline{\scshape{WIGNER'S UNITARY-ANTIUNITARY THEOREM FOR MODULES}}

\pagestyle{myheadings}
\markboth{\textsc{\SMALL WIGNER'S UNITARY-ANTIUNITARY THEOREM FOR
MODULES}}{\textsc{\SMALL WIGNER'S UNITARY-ANTIUNITARY THEOREM FOR
MODULES}}

\normalsize
\vskip .5cm
\centerline{\textsc{Abstract}}

\noindent
Let  $\cH$ be a Hilbert $C^*$-module over a matrix algebra $A$. It is
proved
that any function $T:\cH \to \cH$ which preserves the absolute value
of the (generalized) inner product is of the form $Tf=\varphi(f)Uf$
$(f\in \cH)$, where $\varphi$ is a phase-function and $U$ is an
$A$-linear isometry. The result gives
a natural extension of Wigner's classical unitary-antiunitary theorem
for Hilbert modules.

\vskip .5cm \noindent
1991 \textit{Physics and Astronomy Classification Scheme (PACS)}:
02.30.Sa; 02.30.Tb.

\newpage
\normalsize
\section{Introduction and statement of the result}

Wigner's unitary-antiunitary theorem reads as follows.
Let $H$ be a complex Hilbert space and let $T:H \to H$ be a bijective
function (linearity or continuity is not assumed) with the property that
\[
|\la Tx,Ty\ra |=|\la x,y\ra | \qquad (x,y \in H).
\]
Then $T$ is of the form
\[
Tx=\varphi(x)Ux \qquad (x\in H),
\]
where $U:H \to H$ is either a unitary or an antiunitary operator
and $\varphi :H \to \mathbb C$ is a so-called phase-function which means
that its values are of modulus 1.
This celebrated result plays a very important role in
quantum mechanics and in representation theory in physics.

In our recent paper \cite{MolJAMS} we presented a new, algebraic
approach to this theorem. Our idea turned out to be strong enough
to give a natural generalization of
Wigner's theorem for Hilbert $C^*$-modules over matrix algebras.
However,
in the main result \cite[Theorem 1]{MolJAMS} we supposed that our map is
surjective and, in addition, a  condition was imposed on the underlying
module which was proved to be equivalent to that its so-called modular
dimension is high enough. In the present paper, refining and modifying
our argument quite significantly, we obtain our Wigner-type result in
full generality, that is, neither the surjectivity of the transformation
in question nor the high dimensionality of the Hilbert module is
assumed.

First we clarify the concepts and  notation that we are going to use
throughout. For a bit more detailed discussion we refer to the
introduction of \cite{MolJAMS}. Let $A$ be a $C^*$-algebra.
Let $\cH$ be a left $A$-module with a map $[.,.] : \cH \times \cH \to A$
satisfying
\begin{itemize}
\item[(i)]      $[f+g,h]=[f,h]+[g,h]$
\item[(ii)]     $[af,g]=a[f,g]$
\item[(iii)]    $[g,f]=[f,g]^*$
\item[(iv)]     $[f,f]\geq 0$ and $[f,f]=0$ if and only if $f=0$
\end{itemize}
for every $f,g,h\in \cH$ and $a\in A$.
If $\cH$ is complete with respect to the norm $f \mapsto \|
[f,f]\|^{1/2}$, then we
say that $\cH$ is a Hilbert $A$-module or a Hilbert $C^*$-module over
$A$ with generalized inner product $[.,.]$.
Nowadays, Hilbert modules over $C^*$-algebras
play very important role in many parts of functional analysis such
as, for example, in the K-theory of $C^*$-algebras.
There is another concept of Hilbert modules due to Saworotnow
\cite{Saw}. These are modules over $H^*$-algebras.
The only formal difference in the definition is that in the
case of Saworotnow's modules, the generalized inner product takes its
values in the trace-class of the underlying $H^*$-algebra
and the norm with respect to which we require completeness is
$f \mapsto (\tr [f,f])^{1/2}$.
Saworotnow's modules appear naturally when dealing with
multivariate stochastic processes and
they have applications in Clifford analysis
and hence in some parts of mathematical physics.

If the underlying $C^*$-algebra $A$ is the algebra
$M_d(\mathbb C)$ of all $d\times d$ complex matrices, then, $A$ being
finite dimensional, the norms on $A$ are all equivalent. Therefore,
the Hilbert $C^*$-modules over the $C^*$-algebra $M_d(\mathbb C)$ are
the same as Saworotnow's Hilbert modules over the $H^*$-algebra
$M_d(\mathbb C)$. We emphazise this fact since, in general, the
behaviour
of Saworotnow's Hilbert modules is much nicer and we shall use several
results concerning them.
Finally, we note that it seems to be more common to use right modules
instead of left ones. Of course, this is not a real difference, only a
question of taste.

Now we are in a position to formulate the main result of the paper.
Recall that in any $C^*$-algebra $A$, the element $|a|$ denotes the
square root of $a^*a$ $(a\in A)$.

\begin{theorem}\label{T:wigner}
Let $\cH$ be a Hilbert $C^*$-module over the matrix algebra
$A=M_d(\mathbb C)$, $d>1$.
Let $T:\cH \to \cH$ be a function with the property that
\begin{equation}\label{E:wigner}
|[Tf,Tf']|=|[f,f']| \qquad (f,f' \in \cH).
\end{equation}
Then there exist an $A$-isometry $U:\cH \to \cH$ and a phase-function
$\varphi :\cH \to \mathbb C$ such that
\[
Tf=\varphi(f)Uf \qquad (f\in \cH).
\]
Here, $A$-isometry means that $U:\cH \to \cH$ is a linear map
with $U(af)=aUf$ and $[Uf,Uf']=[f,f']$ $(a\in A,\, f,f'\in \cH)$.
\end{theorem}

The corresponding result for the case $d=1$, that is, when $\cH$ is a
Hilbert space, can be found in \cite{Bar}, \cite{ShA1} (for a recent
paper also see \cite{Rat}).
As we shall see in the proof, the nonappearence of $A$-antiisometries
in the above result is the consequence of the noncommutativity of the
underlying algebra $A$.

Hilbert spaces over algebras different from $\mathbb R$ and
$\mathbb C$ do appear in mathematical physics (see, for example,
\cite{ShA2} for a Wigner-type theorem concerning Hilbert spaces over the
skew-field of quaternions). We believe that our present result may also
have physical interpretation.

\vskip 1cm
\section{Proof}

We give some additional definitions and notation that we shall use in
the proof of our theorem.
As mentioned in the introduction, Saworotnow's modules have
many convenient properties which are familiar in the theory of Hilbert
spaces (we refer to \cite{Saw}). First of all, if $\cH$ is a Hilbert
module over an $H^*$-algebra, then $\cH$ is a Hilbert
space with the inner product $\la .,.\ra =\tr [.,.]$.
If $\cM \subset \cH$ is a closed submodule, then its orthogonal
complement with respect to $\la .,.\ra$ and $[.,.]$ are the same.
A linear operator $T$ on $\cH$ which is bounded with respect
to the Hilbert space norm defined above is called an $A$-linear
operator if $T(af)=aTf$
holds true for every $f\in \cH$ and $a\in A$.
Every $A$-linear operator $T$ is
adjointable, namely, the adjoint $T^*$ of $T$ in the Hilbert space
sense is $A$-linear and we have $[Tf,g]=[f,T^*g]$ $(f,g\in
\cH)$.
Consequently, the collection of all $A$-linear operators forms a
$C^*$-subalgebra
in the full operator algebra on the Hilbert space $\cH$. This will be
denoted by $\cBH$ while the notation of the full operator algebra over a
Hilbert space $H$ is $B(H)$.

In the case of a Hilbert module $\cH$ over an $H^*$-algebra, the
natural equivalent of the Hilbert basis is the so-called modular basis
\cite{MolCZ}. An element $f\in \cH$ is called a modular unit vector, if
$[f,f]$ is a nonzero minimal projection in $A$. A family
$\{f_\alpha\}_\alpha \subset \cH$ is said to be modular orthonormal if

\hskip .5cm  (a) $[f_\alpha, f_\beta ]=0$ if $\alpha\neq \beta$,

\hskip .5cm  (b) $f_\alpha$ is a modular unit vector for every $\alpha$.

\noindent
A maximal modular orthonormal family of vectors in $\cH$ is called a
modular basis.
The common cardinality of modular bases in $\cH$ is called the modular
dimension of $\cH$ (see \cite[Theorem 2]{MolCZ}).

Now, we define operators which are the natural
equivalent of the finite rank operators in the case of Hilbert spaces.
If $f,g \in
\cH$, then let $f\od g$ denote the $A$-linear operator defined by
\[
(f \od g)h=[h,g]f \qquad (h \in \cH).
\]
It is easy to see that for every $A$-linear operator $S$ we
have
\begin{equation*}\label{E:prop1}
S(f\od g)=(Sf)\od g ,\qquad (f\od g) S=f\od (S^*g)
\end{equation*}
and
\begin{equation*}\label{E:prop2}
(f\od g)(f'\od g')=([f',g]f)\od g'=f\od ([g,f']g').
\end{equation*}
Define
\[
\cF(\cH) = \{ \sum_{k=1}^n f_k\od g_k  \, : \, f_k, g_k \in \cH\,
(k=1,\ldots ,n),\,\, n\in \mathbb N \}
\]
which is a *-ideal in the $C^*$-algebra of all $A$-linear operators.
Observe that if $\cH$ is a Hilbert module over $M_d(\mathbb C)$, then
the range of every element of $\cF(\cH)$ has finite linear dimension,
but
there can be finite rank operators on the Hilbert space $\cH$ which do
not belong to $\cF(\cH)$. In general, if the underlying $H^*$-algebra
is infinite dimensional, then these two classes of operators have
nothing to do with each other.

We begin with some auxiliary results that we shall need in the proof of
our theorem.

\begin{lemma}\label{L:lemmus2}
Let $A=M_d(\C)$, $d \in \N$. If $\cH$ is a Hilbert $A$-module, then
every projection in $\cBH$ is of the form $P=\sum_\alpha f_\alpha \od
f_\alpha$, where $\{ f_\alpha \}_\alpha \subset \cH$ is a modular
orthonormal basis in the range of $P$ (the range of an $A$-linear
projection is a closed submodule).

If $\{ f_\alpha \}_\alpha \subset \cH$ is a modular orthonormal set,
then for the orthogonal projection onto the closed submodule generated
by
$\{ f_\alpha \}_\alpha$ (which is an $A$-linear projection) we have
$P=\sum_\alpha f_\alpha \od f_\alpha$.
\end{lemma}

\begin{proof}
Let first $P\in \cBH$ be a projection and let $\{ f\al\}\al$ denote a modular
orthonormal basis in the closed submodule $\rng P$. By
\cite[Theorem 1]{MolCZ} we have
\[
f=\sum_\alpha [f,f\al]f\al \qquad (f \in \rng P)
\]
Since $Pf=0$ and $[f,f\al]=0$ for $f\in \rng P^\perp$, we obtain
$P=\sum_\alpha f_\alpha \od f_\alpha$.

Now, let $\{ f\al \}\al\subset \cH$ be a modular orthonormal set and
denote $\cM$ the closed submodule generated by this set.
We show that
$\{ f\al \}\al$ is a modular basis in $\cM$. Since this collection is
a modular orthonormal family, if this was not maximal, then we
could find a nonzero element $f\in \cM$ which is modular orthogonal to
$\{ f\al \}\al$, that is, $[f,f\al]=0$ for every $\alpha$. But this is a
contradiction, since every element of $\cM$
can be approximated by finite sums of the form $a_1{f\al}_1+\cdots +
a_n{f\al}_n$ $(a_i \in A)$ and hence we would obtain that $f$ is
modular orthogonal to itself. By the first part of the
proof we obtain that the orthogonal projection onto $\cM$ is equal to
$\sum_\alpha f_\alpha \od f_\alpha$, so this operator is an $A$-linear
projection.
\end{proof}

\begin{lemma}\label{L:lemmus1}
Let $A=M_d(\C)$, $d\in \N$ and let $\cH$ be a Hilbert $A$-module.
Suppose that $\cM
\subset \cH$ is a closed submodule and $\{ f\al\}\al$ is a modular
orthonormal system generating $\cM$. Then for every $g,h\in \cM$ we have
\begin{itemize}
\item[(i)] $g=\sum\al [g,f\al]f\al$,
\item[(ii)] $[g,h]=\sum\al [g,f\al][f\al, h]$.
\end{itemize}
Moreover, the vector $k\in \cH$ belongs to $\cM$ if and only if
\[
[k,k]=\sum\al [k,f\al][f\al, k].
\]
\end{lemma}

\begin{proof}
See \cite[Theorem 2]{MolCZ} and its proof.
\end{proof}

\begin{proposition}\label{P:propus}
Let $A=M_d(\C)$, $d \in \N$. If $\cH$ is a Hilbert $A$-module, then
$\cBH$ is a type I von Neumann factor. If the modular dimension of $\cH$
is greater than 2, then $\cBH$ is not ismorphic to $M_2(\C)$.
\end{proposition}

\begin{proof}
It is clear that
$\cBH$ is a von Neumann algebra since it is the commutant of the set
$\{ L_a \, :\, L_af=af\,\, (f\in \cH, \, a\in A)\}$ in the full operator
algebra over $\cH$ as a Hilbert space.
To show that $\cBH$ is a factor, it is sufficient to verify that the
central projections in $\cBH$ are all trivial. Let $P\in \cBH$ be a
nonzero central projection. Let $f$ be a modular unit vector in $\rng
P$. For any $a,b\in A$ we have
\[
P \cdot g \od (af) = g\od (af) \cdot P= g\od P(af) = g\od (af).
\]
This implies that
\[
P(b[f ,f]a^* g)=P((g\od (af))bf)=
(g\od (af))bf=b[f,f]a^* g.
\]
The element $[f,f]$ is a rank-one projection. Hence,
every element of $A$ is the sum of $b[f,f]a^*$-type elements and hence
we obtain that $Pg=g$ for every $g\in \cH$. Thus $P=I$. So, $\cBH$ is a
factor. We next prove that $\cBH$ is type I.
Let $f\in \cH$ be a modular unit vector. Since $[f,f]f=f$ (see
\cite[Lemma 1]{MolCZ}), for any $A\in \cBH$ we compute
\[
f\od f \cdot A \cdot f\od f=([Af,f]f)\od f=
([A([f,f]f),[f,f]f]f\od f)=
\]
\[
([f,f][Af,f][f,f]f)\od f=
\lambda (f\od f)
\]
where $\lambda$ is scalar such that
$[f,f][Af,f][f,f]f= \lambda f$
(the existence of such a scalar follows from the fact that
$[f,f]$ is a rank-one
matrix). This shows that the projection $f\od f$ is abelian. So, every
nonzero central projection in $\cBH$ contains a
nonzero abelian projection which means that $\cBH$ is type I.

Suppose that the modular dimension of $\cH$ is greater than 2.
To see that $\cBH$ is not isomorphic to $M_2(\C)$ it is now enough to
show that the linear dimension of $\cBH$ is greater then 4.
Let $\{ f_1,f_2, f_3\}$ be a modular orthonormal set in $\cH$. Denote
$[f_i,f_i]=e_i$.
If $d\geq 2$, then there are elements $a_i,b_i\in A$ such that $\{
e_ia_i, e_ib_i\}$ is independent for every $i=1,2,3$.
It is easy to check that $\{ (a_if_i)\od f_i, (b_if_i) \od f_i \, : \,
i=1,2,3\}$ is linearly independent. Therefore, the algebraic dimension
of $\cBH$ is at least 6.
If $d=1$, then the statement is trivial.
\end{proof}

Let $H$ be a Hilbert space. Recall that if $x,y\in H$, then $x\ot y$
stands for the operator defined by $(x\ot y)(z)=\la z, y\ra x$ $(z\in
H)$. The ideal of all finite rank operators in $B(H)$ is denoted by
$F(H)$.

\begin{lemma}\label{L:lemmus3}
Let $H$ be a Hilbert space. If $\phi: F(H) \to B(H)$ is a
*-ho\-mo\-morp\-hism which preserves the rank-one projections, then
there is an isometry $U\in B(H)$ such that $\phi$ is of the form
\[
\phi(A)=UAU^* \qquad (A \in F(H)).
\]
Similarly, if $\psi: F(H)\to B(H)$ is a *-antihomomorphism preserving
the rank-one projections, then $\psi$ is of the form
\[
\psi(A)=VA^{tr}V^* \qquad (A \in F(H)),
\]
where $V$ is an isometry and ${}^{tr}$ denotes the transpose with
respect to a fixed orthonormal basis in $H$.
\end{lemma}

\begin{proof}
Let $y,z \in H$ be such that $\la \phi(y\ot y)z,z\ra =1$. Define
\[
Ux=\phi(x \ot y)z \qquad (x\in H).
\]
It is easy to see that $U$ is an isometry and $UA=\phi(A)U$ $(A\in
F(H))$. Let $x\in H$ be an arbitrary unit vector. Then $\phi(x\ot x)$ is
a rank-one projection, so it is of the form $\phi(x\ot x)=x'\ot x'$
with some unit vector $x'\in H$. Since
\[
Ux\ot x=\phi(x\ot x)U=x'\ot U^*x',
\]
we obtain that $x'$ is equal to $Ux$ multiplied by a scalar of modulus
1. Therefore, $\phi(x\ot x)=Ux \ot Ux =U\cdot x\ot x \cdot U^*$. Since
this holds true for every unit vector $x\in H$, by linearity we have the
first assertion of the lemma.

As for the second statement, we can apply a similar argument. Choosing
$y,z\in H$ such that $\la \psi(y\ot y)z, z \ra=1$,
define
\[
\tilde V x=\psi(y \ot x)z \qquad (x\in H).
\]
One can verify that $\tilde V$ is an antiisometry (that
is, a conjugate-linear isometry), and
then prove that $\psi(A)=\tilde V A^* {\tilde V}^*$ $(A\in F(H))$.
Considering an antiunitary operator $J$ for which $JA^*J^*=A^{tr}$ and
defining $V=\tilde VJ$, we conclude the proof.
\end{proof}

\begin{lemma}\label{L:lemmus4}
Let $(a_n)$ be a sequence in the Hilbert space $H$ and let $b\in H$ be
such that $\sum_n a_n \ot a_n =b\ot b$ in the trace norm. Then for every
$n$ there exists a scalar $\lambda_n$ such that $a_n =\lambda_n b$.
\end{lemma}

\begin{proof}
Clearly, we may assume that $\| b\|=1$. Taking traces on both sides
of the equality $\sum_n a_n \ot a_n =b\ot b$, we obtain $\sum_n
\| a_n \|^2=1$. On the other hand, we also have
\[
\sum_n |\la b, a_n\ra|^2=
\la(\sum_n a_n \ot a_n)b,b\ra =1.
\]
By Schwarz inequality
\[
1= \sum_n |\la b, a_n\ra|^2\leq \sum_n \| a_n \|^2 =1.
\]
So, there are equalities in the Schwarz inequalites $|\la b,a_n \ra
|\leq \| a_n\|$. This implies the assertion.
\end{proof}

\begin{proof}
We define an orthoadditive projection-valued measure $\mu$ on the
lattice $\cPH$ of all $A$-linear projections
as follows. If $\{ f\al\}\al$ is a modular orthonormal set, then let
\[
\mu(\sum\al f\al \od f\al)=\sum\al Tf\al \od Tf\al.
\]
Observe that by \eqref{E:wigner}, $\{ Tf\al \}\al$ is also modular
orthonormal and, hence, by Lemma~\ref{L:lemmus2} $\sum\al Tf\al \od
Tf\al$ belongs to $\cPH$. We show that $\mu$ is well-defined. Let $\{
f\al \}\al$ and $\{
g\be\}\be$ generate the same closed submodule $\cM$. We claim that the
same holds true for $\{ Tf\al \}\al$ and $\{ Tg\be\}\be$.
Indeed, if $g\in \cM$, then due to the fact that $\{ f\al \}\al$ is a
modular basis in $\cM$ we see that $g=\sum\al [g, f\al ]f\al$. This
implies that
\[
[Tg,Tg]=[g,g]=\sum\al [g,f\al][f\al, g]=
\sum\al [Tg,Tf\al][Tf\al, Tg]
\]
which, by Lemma~\ref{L:lemmus1}, gives us that $Tg$ belongs to the
closed submodule generated by $\{ Tf\al \}\al$.
It is now obvious that $\mu$ is an orthoadditive
$\cPH$-valued measure on $\cPH$.

Let us suppose that the modular dimension of $\cH$ is greater than 2. By
Proposition~\ref{P:propus} we can apply a deep result of Bunce and
Wright
\cite[Theorem A]{BW}. It states that every bounded finitely
orthoadditive, Banach space valued
measure on the set of all projections in a von Neumann algebra without a
summand isomorphic to $M_2(\C)$ can be uniquely extended to a bounded
linear transformation defined on the whole algebra. Let $\phi: \cBH \to
\cBH$ denote the transformation corresponding to $\mu$. Since it sends
projections to projections,
it is a standard argument to verify that $\phi$ is a Jordan
*-endomorphism of $\cBH$, that is, we have $\phi(T)^2=\phi(T^2)$,
$\phi(T)^*=\phi(T^*)$ $(T\in \cBH)$ (see, for example, the proof of
\cite[Theorem 2]{MolStud1}).

We prove that $\phi(f\od f)=Tf \od Tf$ for every $f\in \cH$.
Let $[f,f]=\sum_i \lam_i^2 e_i$, where $\lam_i$'s are nonnegative real
numbers and $e_i$'s are pairwise orthogonal rank-one projections.
Define $f_i=(1/\lam_i)e_if$. We have $[f_i,f_i]=e_i$ and $[f_i,f_j]=0$
if $i\neq j$, that is, $\{ f_i\}_i$ is modular orthonormal. Then
$f=\sum_i \lam_if_i=\sum_i e_i f$ since
\[
\sum_i [f,f_i][f_i,f]=\sum_i \lam_i^2e_i=[f,f]
\]
implies that
$f=\sum_i [f,f_i]f_i=\sum_i e_if$ (see Lemma~\ref{L:lemmus1}).
So, we have
\[
\phi(f\od f)=\sum_{i,j} \phi(e_if \od e_jf).
\]
But $(e_if)\od (e_jf)=0$ if $i\neq j$. Indeed, we compute
$[g,e_jf]e_if=[g,f]e_je_if=0$ for every $g\in \cH$. Hence,
\[
\phi(f\od f)=\sum_{i} \phi(e_if \od e_if)=
\sum_{i} \lam_i^2\phi(f_i \od f_i)=
\]
\[
\sum_{i} \lam_i^2\mu(f_i \od f_i)=
\sum_{i} \lam_i^2 Tf_i \od Tf_i .
\]
So, the question is that whether the equality
$Tf\od Tf=\sum_{i} \lam_i^2 Tf_i \od Tf_i$ holds true.
Clearly, $\{ Tf_i\}$ is modular orthonormal.
We compute
\[
[Tf,Tf]=[f,f]=\sum_i [f,f_i][f_i,f]=
\sum_i [Tf,Tf_i][Tf_i,Tf]
\]
which, by Lemma~\ref{L:lemmus1}, implies that $Tf=\sum_i [Tf,
Tf_i]Tf_i$. We know that $|[Tf,Tf_i]|=|[f,f_i]|=\lam_i e_i$.
Similarly, $|[Tf_i,Tf]|=|[f_i,f]|=\lam_i e_i$.
Since $e_i$ is a rank-one projection, we obtain that $[Tf,Tf_i]$ is also
rank-one. Furthermore, as $|[Tf,Tf_i]|=|[Tf_i,Tf]|$ is a scalar multiple
of $e_i$
we can infer that $[Tf,Tf_i]=\mu_i \lam_i e_i$, where $\mu_i$ is a
scalar of modulus 1.
Therefore, we have
\[
Tf\od Tf=\sum_{i,j} \mu_i {\overline \mu_j} (\lam_i e_i Tf_i \od
\lam_j e_j Tf_j).
\]
But similarly as above, for $i\neq j$ we have
\[
(e_i Tf_i \od e_j Tf_j)g=
[g, e_jTf_j]e_iTf_i=[g, Tf_j]e_je_iTf_i=0.
\]
Therefore
\[
Tf\od Tf=
\sum_{i,j} \mu_i {\overline \mu_j} (\lam_i e_i Tf_i \od \lam_j e_j
Tf_j)=
\sum_{i} \mu_i {\overline \mu_i} (\lam_i e_i Tf_i \od \lam_i e_i Tf_i)=
\]
\[
\sum_{i} \lam_i e_i Tf_i \od \lam_i e_i Tf_i=
\sum_{i} \lam_i^2 (e_i Tf_i \od e_i Tf_i)   .
\]
But
$(e_i Tf_i \od e_i Tf_i)=Tf_i\od Tf_i$. Indeed, since $Tf_i$ is a
modular unit vector, we have
$e_iTf_i=[f_i,f_i]Tf_i=[Tf_i,Tf_i]Tf_i=Tf_i$ (see \cite[Lemma
1]{MolCZ}). Consequently, we obtain
$Tf\od Tf=\sum_{i} \lam_i^2 Tf_i \od Tf_i$ and this was to be proved.
So, we get $\phi(f\od f)=Tf \od Tf$ for every $f\in \cH$.

We assert that $\phi$ is either a *-homomorphism or a
*-anti\-ho\-mo\-mor\-phism. By Lemma~\ref{L:lemmus2} the minimal
projections in $\cH$
are exactly the operators of the form $f\od f$, where $f\in \cH$ is a
modular unit vector. Clearly, $\phi$ sends minimal projections to
minimal
projections. By \cite[Lemma 2]{MolJAMS} the linear space generated by
the minimal projections in $\cBH$ is $\cF(\cH)$.
Since $\cBH$ is a type I factor,
it is isomorphic to the full operator algebra $B(H)$ on a Hilbert space
$H$. Since *-isomorphisms preserve the minimal projections,
$\cF(\cH)$ corresponds to the ideal $F(H)$
of all finite rank operators in $B(H)$. Under this identification,
we obtain a Jordan *-homomorphism $\tilde \phi$ on $F(H)$
corresponding to $\phi_{|\cF(\cH)}$ which sends rank-one projections to
rank-one projections. Since $F(H)$ is a local matrix algebra, by
\cite[Theorem 8]{JR} we obtain that $\tilde \phi$ is the sum of a
*-homomorphism
and a *-antihomomorphism. As $\tilde \phi$ preserves the rank-one
projections, from the simplicity of the ring $F(H)$ it follows that
$\tilde \phi$ is either a *-homomorphism
or a *-antihomomorphism. Obviously, the same holds for
$\phi_{|\cF(\cH)}$.

Let us suppose that the
modular dimension of $\cH$ is greater than $d$. By \cite[Remark
2]{MolJAMS}, there are vectors $g,h\in \cH$ such that $[g,h]=I$.
The map $\phi_{|\cFH}$ is either a *-homomorphism or a
*-antihomomorphism. First consider this latter case.
Referring to Lemma~\ref{L:lemmus3} we have an
operator $U\in \cBH$ with $U^*U=I$ and a *-antiautomorphism $\psi$ of
$\cFH$ such that $\phi(A)=U\psi(A)U^*$ $(A\in \cFH)$.

We define
\[
Vf=\psi(g\od f)U^*Th \qquad (f \in \cH)
\]
where $g,h\in \cH$ are fixed and such that $[g,h]=I$. Clearly, $V$ is a
conjugate-linear operator. We have
\[
VAf=\psi(g\od (Af))U^*Th=
\psi(g\od f A^* )U^*Th=
\]
\[
\psi(A)^*\psi(g\od f)U^*Th=
\psi(A)^*Vf,
\]
that is, $VA=\psi(A)^*V$ $(A\in \cFH)$.
We compute
\[
[Vf,Vf]=
[\psi(g\od f)U^*Th, \psi(g\od f)U^*Th]=
\]
\[
[\psi(g\od f \cdot f\od g)U^*Th, U^*Th]=
[U\psi(g\od f \cdot f\od g)U^*Th, Th]=
\]
\[
[\phi(g\od f \cdot f\od g)Th, Th]=
[\phi(\sqrt{[f,f]}g\od \sqrt{[f,f]}g)Th, Th]=
\]
\[
[(T(\sqrt{[f,f]}g)\od T(\sqrt{[f,f]}g))Th, Th]=
\]
\[
[Th, T(\sqrt{[f,f]}g)][T(\sqrt{[f,f]}g), Th]=
\]
\[
[h, \sqrt{[f,f]}g][\sqrt{[f,f]}g, h]=
[h, g][f,f][g, h]=
[f,f].
\]
Since $V$ is conjugate-linear, by polarization we obtain
\[
[Vf,Vf']=[f',f] \qquad (f, f'\in \cH).
\]

We show that $\rng T\subset \rng U$ which will imply $UU^*T=T$ ($UU^*$
is the projection onto the range of $U$). Let $f\in
\cH$. In the previous part of the proof we have learnt that $Tf
\od Tf$ is a linear combination
of operators of the form $Tf_b \od Tf_b$, where $f_b$'s are modular unit
vectors. We have
\[
Tf_b \od Tf_b=\phi(f_b \od f_b)=U\psi(f_b \od f_b) U^*
\]
and, $\psi$ being a *-antiautomorphism, $\psi(f_b\od f_b)$ is a minimal
projection. Therefore,
$\psi(f_b\od f_b)=f'_b \od f'_b$ with some modular unit vector
$f'_b$ and hence $Tf_b\od Tf_b=Uf'_b \od Uf'_b$. Now let
$Tf=g'+g''$, where
$g' \in \rng U$ and $g''\in \rng U^\perp$. We have
\[
[g'',g'']^2=[g'',Tf][Tf,g'']=[(Tf\od Tf)g'',g'']=0.
\]
This gives us that $g''=0$ which shows that $Tf \in \rng U$.

We next prove that $V$ is surjective.
Let $f\in \cH$ be arbitrary. Since $\psi$ is a *-antiautomorphism of
$\cFH$, we
can find an operator $R\in \cFH$ such that $\psi(R)^*=f\od U^*Th$.
We compute
\[
VRg=
\psi(R)^*Vg=
\psi(R)^*\psi(g\od g)U^*Th=
\]
\[
\psi(R)^*U^*\phi(g\od g)Th=
\psi(R)^*U^*(Tg\od Tg)Th=
[Th, Tg]\psi(R)^*U^*Tg=
\]
\[
[Th, Tg][U^*Tg, U^*Th]f=
[Th, Tg][UU^*Tg, Th]f=
\]
\[
[Th, Tg][Tg, Th]f=
[h, g][g, h]f=f.
\]
Since $f$ was arbitrary, we have the surjectivity of $V$.

We compute
\[
[UVf', Tf][Tf, UVf']=
[(Tf\od Tf)UVf', UVf']=
\]
\[
[U^*(Tf\od Tf)UVf', Vf']=
[U^*\phi(f\od f)UVf', Vf']=
\]
\[
[\psi(f\od f)Vf', Vf']=
[(V\cdot f\od f)f', Vf']=
\]
\[
[f', (f\od f)f']=
[f',f][f,f'].
\]
This gives us that
\[
[V^{-1}U^*Tf',f][f,V^{-1}U^*Tf']=
[UVV^{-1}U^*Tf',Tf][Tf,UVV^{-1}U^*Tf']=
\]
\[
[UU^*Tf',Tf][Tf,UU^*Tf']=
[Tf',Tf][Tf,Tf']=
[f',f][f,f'].
\]
Replacing $f$ by $xf$ $(x\in A)$, we obtain
\[
[V^{-1}U^*Tf',f]x^*x[f,V^{-1}U^*Tf']=
[f',f]x^*x[f,f'].
\]
Since every element of $A$ is a linear combination of elements of the
form $x^*x$, it follows that
\[
[V^{-1}U^*Tf',f]y[f,V^{-1}U^*Tf']=
[f',f]y[f,f']
\]
holds for every $y\in A$. This implies that for every $f\in \cH$, the
matrices $[f,V^{-1}U^*Tf']$ and $[f,f']$ are linearly dependent.
It requires only elementary linear algebra to verify the following
assertion. If $X,Y$ are vector spaces and $A,B :X \to Y$ are linear
operators such that for every $x\in X$, the set $\{ Ax, Bx\}$ is
linearly
dependent, then either $A$ and $B$ have rank at most one or $\{
A,B\}$ is linearly dependent.
Since the rank
of the linear operator $f\mapsto [f,f']$ is clearly greater than 1 if
$f'\neq 0$, we have a scalar $\lambda_{f'}$ (depending only on $f'$)
such that $[f,V^{-1}U^*Tf']=\lambda_{f'}[f,f']$ $(f,f'\in \cH)$.
This gives us that there is a function $\varphi :\cH \to \C$ such that
$V^{-1}U^*Tf'={\overline{\varphi(f')}}f'$ which results in
$Tf'=\varphi(f')UVf'$. It follows from the properties
of $T,U,V$ that $\varphi$ is of modulus 1.
Finally, we have
\[
|[f,f']|=|[Tf,Tf']|=|[UVf,UVf']|=
|[Vf,Vf']|=|[f',f]|.
\]
Since this must hold true for every $f,f'\in \cH$, it follows that for
every rank-one matrix $a\in A$ we have $|a|=|a^*|$. But this is an
obvious contradiction. Since we have started with assuming that
$\phi_{|\cFH}$ is a *-antihomomorphism, we thus obtain that
it is in fact a *-homomorphism.

Pushing the problem from $\cBH$ to the full operator algebra $B(H)(\cong
\cBH)$, we see that there is an
$A$-isometry $U\in \cBH$ such that $\phi(A)=UAU^*$ $(A\in \cFH)$. This
gives us
that $Tf\od Tf=Uf\od Uf$ for every $f\in \cH$. Similarly as before, this
implies that $\rng T\subset \rng U$ which yields
$UU^*Tf=Tf$ $(f\in \cH)$. We next compute
\[
[Uf', Tf][Tf, Uf']=
[(Tf\od Tf)Uf',Uf']=
\]
\[
[(Uf\od Uf)Uf',Uf']=
[Uf', Uf][Uf, Uf']=
[f', f][f, f'],
\]
which gives us that
\[
[U^*Tf',f][f,U^*Tf']=[UU^*Tf',Tf][Tf,UU^*Tf']=
\]
\[
[Tf',Tf][Tf,Tf']=[f',f][f,f'].
\]
Just as above, it follows that $U^*Tf'$ is a scalar multiple of $f'$.
Therefore,
there exists an $A$-isometry $U$ and a phase-function $\varphi :\cH \to
\C$ such that
\[
Tf =\varphi(f)Uf \qquad (f\in \cH).
\]
This completes the proof in the case when the modular dimension $n$ of
$\cH$ is greater than $d$.

We now treat the low dimensional cases, that is, when
$n\leq d$.
Let $H_d$ denote the $d$-dimensional complex Euclidean space. Then $H_d$
can be considered as a Hilbert $A$-module. Here, the module
operation is $(a,\xi)\mapsto a(\xi)$ and the generalized inner product
is defined by $[\xi,\zeta]=\xi \ot \zeta$.
Clearly, the modular dimension of this module is 1. It now follows
from the structure of our Hilbert $A$-modules (see, for example,
\cite{Oza}) that $\cH$ is isomorphic to the
$n$-fold direct sum of $H_d$ with itself. So, we may assume that $\cH
=\sum_{i=1}^n \oplus H_d$. The definition of
the module operation and that of the inner product on this direct sum is
defined as follows
\[
a[\xi_i]_i=[a\xi_i]_i, \qquad
[ [\xi_i]_i, [\zeta_i]_i]=
\sum_i \xi_i \ot \zeta_i.
\]
Let us describe the elements of $\cBH$. Since every element of $\cBH$
is a linear operator on the direct sum of vector spaces, it can
represented by a matrix
\[
\left[
\begin{matrix}
a_{11} &\dots& a_{1n}\\
\vdots &\ddots&\vdots \\
a_{n1} &\dots& a_{nn}
\end{matrix}
\right]
\]
where $a_{ij}$'s are linear operators acting on $H_d$. Now,
$A$-linearity means that
\[
\left[
\begin{matrix}
a_{11}a\xi_1 +\cdots + a_{1n}a\xi_n\\
\vdots  \\
a_{n1}a\xi_1 +\cdots + a_{nn}a\xi_n
\end{matrix}
\right]
=
\left[
\begin{matrix}
a_{11} &\dots& a_{1n}\\
\vdots &\ddots&\vdots \\
a_{n1} &\dots& a_{nn}
\end{matrix}
\right]
\left[
\begin{matrix}
a\xi_1\\
\vdots\\
a\xi_n
\end{matrix}
\right]
=
\]
\[
\left[
\begin{matrix}
a(a_{11}\xi_1 +\cdots + a_{1n}\xi_n)\\
\vdots  \\
a(a_{n1}\xi_1 +\cdots + a_{nn}\xi_n)
\end{matrix}
\right]
\]
holds for every $a\in A$ and $\xi_i\in H_d$. It is easy to see that this
is equivalent to $a_{ij}a=aa_{ij}$ $(a\in A)$ which means that
$a_{ij}$'s are scalars. Consequently, $\cBH$ is isomorphic
to $M_n(\C)$.

Suppose that $n>1$.
If $\zeta$ is any vector in $H_d$, then let $\zeta^k$ denote the element
of $\cH$ whose coordinates are all 0 except for the $k$th one which is
$\zeta$. Fix a unit vector $\xi \in H_d$.
We have
\[
\sum_i (T\xi^k)_i \ot (T\xi^k)_i=
[T\xi^k ,T\xi^k]=[\xi^k, \xi^k]=\xi \ot \xi.
\]
From Lemma~\ref{L:lemmus4}
we infer that for every $i=1, \ldots, n$, there is a scalar
$\alpha_{ik}$ such that $(T\xi^k)_i=\alpha_{ik}\xi$.
Clearly, the columns of the matrix $(\alpha_{ik})$ are unit
vectors.
Since $[T\xi^k, T\xi^l]=0$ for $k\neq l$, it follows that the columns of
our matrix are pairwise orthogonal as well. So $(\alpha_{ik})$ is a
unitary matrix and hence it defines an $A$-unitary operator $U$ on
$\cH$.
Considering $U^*T$ instead of $T$, we can assume that $T\xi^k$ is equal
to $\xi^k$ for every $k=1, \ldots ,n$.
If $f$ is any vector in $\cH$, then considering the equality
\[
|\xi\ot (Tf)_k|=|[T\xi^k, Tf]|=|[\xi^k, f]|=|\xi \ot f_k|
\]
we obtain
\begin{equation}\label{E:reka1}
(Tf)_k=\mu_k f_k \qquad (k=1, \ldots ,n)
\end{equation}
with some
scalars $\mu_k$ of modulus 1. We claim that all the $\mu_k$'s are equal.
Fix a $g\in \cH$ whose coordinates are pairwise orthogonal unit
vectors in $H_d$ (recall that $n\leq  d$).
It is apparent that if we multiply $T$ from the left by an $A$-unitary
operator whose matrix is diagonal, then the so obtained transformation
still has the property \eqref{E:reka1}. So we may assume that $Tg=g$.
Let $f\in \cH$ be arbitrary. We have
\[ | \sum_i \mu_i f_i\ot g_i|=
|[Tf , Tg]|=|[f , g]|=
| \sum_i f_i\ot g_i|.
\]
This implies that
\[
\sum_{i,j} \la \mu_j f_j, \mu_i f_i\ra g_i \ot g_j=
\sum_{i,j} \la f_j, f_i\ra g_i \ot g_j
\]
which gives that
\[
\la \mu_j f_j, \mu_i f_i\ra=
\la f_j, f_i\ra.
\]
So, if $\la f_i, f_j\ra\neq 0$, then we have $\mu_i=\mu_j$. Suppose
now that $\la f_i,f_j \ra=0$ but $f_i,f_j \neq 0$.
Let $\zeta \in H_d$ be any nonzero vector and consider
$\zeta^i+\zeta^j$. By what
we have just proved, it follows that  $T(\zeta^i+\zeta^j)$ is a scalar
multiple of $\zeta^i+\zeta^j$. We compute
\[
|\zeta\ot (\mu_if_i + \mu_jf_ j)|=
|[\zeta^i+\zeta^j, Tf]|=
|[T(\zeta^i+\zeta^j), Tf]|=
\]
\[
|[\zeta^i+\zeta^j, f]|=
|\zeta\ot (f_i + f_j)|
\]
which clearly gives us that $\mu_i=\mu_j$. Therefore, we obtain
that for any vector $f\in \cH$, $Tf$ is equal to $f$ multiplied by a
complex number of modulus 1. The assertion of the theorem now follows
for the case $1<n \leq d$.

Finally, suppose that $n=1$, which means that $\cH=H_d$. Our problem is
to describe those maps
$T:H_d \to H_d$ for which $|T\xi \ot T\zeta|=|\xi \ot \zeta|$ $(\xi,
\zeta \in H_d)$. But this equality clearly implies that $T\zeta$ is
equal to $\zeta$ multiplied by a scalar of modulus 1.

The proof of the theorem is now complete.
\end{proof}

\section{Acknowledgements}
This research was supported from the following sources:
(1) Hungarian National Foundation for Scientific Research
(OTKA), Grant No. T--030082 F--019322,
(2) a grant from the Ministry of Education, Hungary, Reg.
No. FKFP 0304/1997.

\newpage

\end{document}